\input amstex
\input amsppt.sty
\magnification=\magstep1
\TagsOnRight
\pageheight{45pc}
\nologo
\def\N{\Bbb N}
\def\Z{\Bbb Z}

\def\l{\left}
\def\r{\right}
\def\bg{\bigg}
\def\({\bg(}
\def\[{\bg[}
\def\){\bg)}
\def\]{\bg]}
\def\t{\text}
\def\f{\frac}
\def\mo{\roman{mod}}
\def\em{\emptyset}

\def\eq{\equiv}

\def\ls{\leqslant}
\def\gs{\geqslant}
\def\al{\alpha}

\def\Proof{\noindent{\it Proof}}

\def\Remark{\medskip\noindent{\it  Remark}}
\def\Ack{\medskip\noindent {\bf Acknowledgment}}

\topmatter \hbox{Math. Comp., in press.}
\bigskip
\title  Covers of the integers with odd moduli and their applications
to the forms $x^m-2^n$ and $x^2-F_{3n}/2$\endtitle
\rightheadtext{Covers of the integers with odd moduli and their applications}
\author Ke-Jian Wu$^1$  and Zhi-Wei Sun$^{2,*}$\endauthor
\leftheadtext{Ke-Jian Wu and Zhi-Wei Sun}
\affil $^1$Department of
Mathematics, Zhanjiang Normal University
\\Zhanjiang 524048, People's Republic of China
\\ {\tt kjwu328\@yahoo.com.cn}
\medskip
$^2$Department of Mathematics, Nanjing University
\\Nanjing 210093, People's Republic of China
\\ {\tt zwsun\@nju.edu.cn}
\\ http://math.nju.edu.cn/$\sim$zwsun\endaffil

\abstract In this paper we construct a cover $\{a_s(\mo\
n_s)\}_{s=1}^k$ of $\Z$ with odd moduli such that there are
distinct primes $p_1,\ldots,p_k$ dividing
$2^{n_1}-1,\ldots,2^{n_k}-1$ respectively. Using this cover we
show that for any positive integer $m$ divisible by none of
$3,\, 5,\, 7,\,11,\,13$
there exists an infinite arithmetic progression of positive odd
integers the $m$th powers of whose terms are never of the form
$2^n\pm p^a$ with $a,n\in\{0,1,2,\ldots\}$ and $p$ a prime. We
also construct another cover of $\Z$ with odd moduli and use it to
prove that $x^2-F_{3n}/2$ has at least two distinct prime factors
whenever $n\in\{0,1,2,\ldots\}$ and $x\eq a\ (\mo\ M)$, where
$\{F_i\}_{i\gs0}$ is the Fibonacci sequence, and $a$ and $M$ are
suitable positive integers having 80 decimal digits.
\endabstract
\thanks {\it Keywords}: Cover of the integers, arithmetic progression, Fibonacci sequence,
prime divisor.
\newline\indent
2000 {\it Mathematics Subject Classification}.
Primary 11B25; Secondary 11A07, 11A41, 11B39, 11D61, 11Y99.
\newline\indent *This author is responsible for communications,
and supported by the National Natural Science Foundation (grant
10871087) of China.
\endthanks
\endtopmatter

\document

\heading {1. Introduction}\endheading

For $a\in\Z$ and $n\in\Z^+=\{1,2,3,\ldots\}$ we let
$$a(n)=\{x\in\Z:\,x\equiv a\,(\mo\,n)\}$$
which is a residue class modulo $n$. A finite system
$$A=\{a_s(n_s)\}^k_{s=1}\tag1.1$$
of residue classes is said to be a {\it cover} of $\Z$ if every
integer belongs to some members of $A$. Obviously (1.1) covers all
the integers if it covers $0,1,\ldots,N_A-1$ where
$N_A=[n_1,\ldots,n_k]$ is the least common multiple of the moduli
$n_1,\ldots,n_k$. The reader is referred to [Gu] for problems and
results on covers of $\Z$ and to [FFKPY] for a recent breakthrough
in the field. In this paper we are only interested in applications
of covers.

 By a known result of Bang [B] (see also Zsigmondy [Z]
 and Birkhoff and Vandiver [BV]), for each integer $n>1$ with
 $n\not=6$, there exists a prime factor
 of $2^n-1$ not dividing $2^m-1$ for any $0<m<n$; such a prime is
 called a {\it primitive prime divisor of $2^n-1$}. P. Erd\H os,
who introduced covers of $\Z$ in the early 1930s, constructed the
following cover (cf. [E])
$$A_0=\{0(2),\ 0(3),\ 1(4),\ 3(8),\ 7(12),\ 23(24)\}$$
whose moduli are distinct, greater than one and different from 6.
It is easy to check that
$2^2-1,2^3-1,2^4-1,2^8-1,2^{12}-1,2^{24}-1$ have primitive prime
divisors $3,7,5,17,13,241$ respectively. Using the cover $A_0$ and
the Chinese Remainder Theorem, Erd\H os showed that any integer
$x$ satisfying the congruences
$$\cases x\eq2^0\ (\mo\ 3),&\\x\eq2^0\ (\mo\ 7),&\\x\eq2^1\ (\mo\ 5),&
\\x\eq2^3\ (\mo\ 17),&\\x\eq 2^7\ (\mo\ 13),&\\x\eq2^{23}\ (\mo\ 241)&\endcases$$
and the additional congruences $x\eq1\ (\mo\ 2)$ and $x\eq 3\
(\mo\ 31)$ cannot be written in the form $2^n+p$ with
$n\in\N=\{0,1,2,\ldots\}$ and $p$ a prime. The reader may consult
[SY] for a refinement of this result. By improving the work of
Cohen and Selfridge [CS], Sun [S00] showed that for any integer
$$x\eq 47867742232066880047611079\ \(\mo\ \prod_{p\in P}p\)$$
with
$$P=\{2,\,3,\,5,\,7,\,11,\,13,\,17,\,19,\,31,\,37,\,41,\,61,\,73,\,97,\,109,\,151,\,241,\,257,\,331\},$$
we have $x\not=\pm p^a\pm q^b$ where $p,q$ are primes and
$a,b\in\N$. In 2005, Luca and St\u{a}nic\u{a} [LS] constructed a
cover of $\Z$ to show that if $n$ is sufficiently large and $n\eq
1807873\ (\mo\ 3543120)$ then $F_n\ne p^a+q^b$ with $p, q$ prime
numbers and $a, b\in\N$, where the Fibonacci sequence
$\{F_n\}_{n\gs0}$ is given by
$$F_0=0,\ F_1=1,\ \t{and}\ F_{n+1}=F_n+F_{n-1}\ \t{for}\ n=1,2,3,\ldots.$$

A famous conjecture of Erd\H os and J. L. Selfridge states that
there does not exist a cover of $\Z$ with all the moduli odd,
distinct and greater than one. There is little progress on this
open conjecture (cf. [Gu] and [GS]). In contrast, we have the
following theorem.

\proclaim{Theorem 1.1} There exists a cover $A_1=\{a_s(n_s)\}_{s=1}^{173}$ of $\Z$
with all the moduli greater than one and dividing the odd number
$$3^3\times 5^2\times7\times11\times 13=675675,$$
for which there are distinct primes $p_1,\ldots,p_{173}$ greater
than $5$ such that each $p_s\ (1\ls s\ls 173)$ is a primitive
prime divisor of $2^{n_s}-1$.
\endproclaim

Theorem 1.1 has the following application.

\proclaim{Theorem 1.2} Let $N$ be any positive integer. Then there
is a residue class consisting of odd numbers such that for each
nonnegative $x$ in the residue class and each $m\in\{1,\ldots,N\}$
divisible by none of $3,\,5,\,7,\,11,\,13$, the number $x^m-2^n$
with $n\in\N$ always has at least two distinct prime factors.
\endproclaim

\Remark\ 1.1. Let $m\in\Z^+$. Chen [C] conjectured that there are
infinitely many positive odd numbers $x$ such that $x^m-2^n$ with
$n\in\Z^+$ always has at least two distinct prime factors, and he
was able to prove this when $m\eq1\ (\mo\ 2)$ or $m\eq\pm 2\ (\mo\
12)$. The conjecture is particularly difficult when $m$ is a high power of 2.
In a recent preprint [FFK], Filaseta, Finch and Kozek
confirmed the conjecture for $m=4,6$ with the help of
a deep result of Darmon and Granville [DG] on generalized Fermat
equations; they also showed that there exist
infinitely many integers $x\in\{1,3^8,5^8,\ldots\}$ such that
$x^m2^n+1$ with $n\in\Z^+$ always has at least two distinct prime
divisors.
\medskip

Recall that $\{F_n\}_{n\gs0}$ is the Fibonacci sequence. Set
$u_n=F_{3n}/2$ for $n\in\N$. Clearly, $u_0=0$, $u_1=1$, and
$$\align u_{n+1}=&\f{F_{3n+3}}2=\f{F_{3n+1}+(F_{3n+1}+F_{3n})}2
\\=&F_{3n+1}+u_n=F_{3n-1}+3u_n
\\=&4u_n+\f{2F_{3n-1}-F_{3n}}2
\\=&4u_n+\f{F_{3n-1}-F_{3n-2}}2
\\=&4u_n+u_{n-1}
\endalign$$
for every $n=1,2,3,\ldots$.

Now we give the third theorem which is of a new type
and will be proved on the basis of certain cover of $\Z$ with odd moduli.

\proclaim{Theorem 1.3} Let
$$\align
a=&312073868852745021881735221320236651673651936708237682^{\_}\\
&34185354856354918873864275
\endalign$$
and
$$\align M=&368128524439220711844024989130760705031462298208612115^{\_}\\
&58347078871354783744850778.
\endalign$$
Then, for any $x\eq a\ (\mo\ M)$ and $n\in\N$, the number
$x^2-F_{3n}/2$ has at least two distinct prime divisors.
\endproclaim

\Remark \ 1.2. (a) Actually our proof of Theorem 1.3 yields the
following stronger result: Whenever $y\in a^2(M)$ and $n\in\N$,
the number $y-F_{3n}/2$ has at least two distinct prime divisors.

(b) In view of Theorem 1.3, it is interesting to study the
diophantine equation $x^2-F_{3n}/2=\pm p^a$ with $a,n,x\in\N$ and
$p$ a prime, or the equation $F_{3n}=2x^2\pm dy^2$ with $d$ equal
to $1$ or $2$ or twice an odd prime. The related equation
$F_n=x^2+dy^2$ has been investigated by Ballot and Luca [BL].

\bigskip

The second author has the following conjecture.

\proclaim{Conjecture 1.1} Let $m$ be any positive integer. Then
there exist $b,d\in\Z^+$ such that whenever $x\in b^m(d)$ and
$n\in\N$ the number $x-F_n$ has at least two distinct prime
divisors. Also, there are odd integer $b$ and even number $d\in\Z^+$
such that whenever $x\in b^m(d)$ and
$n\in\N$ the number $x-2^n$ has at least two distinct prime
divisors.
\endproclaim

\Remark\ 1.3. (a) We are unable to prove Conjecture 1.1 since it
is difficult for us to construct a suitable cover of $\Z$ for the
purpose.

(b) In 2006, Bugeaud, Mignotte and Siksek [BMS] showed that the
only powers in the Fibonacci sequence are
$$F_0=0,\ F_1=F_2=1,\ F_6=2^3\ \t{and}\ F_{12}=12^2.$$
It seems challenging to solve the diophantine equation
$x^m-F_n=\pm p^a$ with $a,n,x\in\N$, $m>1$, and $p$ a prime.

\medskip

We are going to show Theorems 1.1--1.3 in Sections 2--4 respectively.

\heading{2. Proving Theorem 1.1 via constructions}\endheading

\noindent{\it Proof of Theorem 1.1}. Let $a_1(n_1),\ldots,a_{173}(n_{173})$
be the following 173
residue classes respectively.

$$\aligned
&0(3),\,1(5),\,0(7),\,1(9),\,7(11),\,8(11),\,7(13),\,8(15),\,
19(21),\,17(25),\,22(25),\\
&25(27),\,23(33),\,29(35),\,30(35),\,
14(39),\,17(39),\,4(45),\,13(45),\,0(55),\\
&25(55),\,50(55),\,
25(63),\,52(63),\,9(65),\,2(75),\,32(75),\,13(77),\,41(91),\\
&62(91),\,76(91),\,5(99),\,65(99),\,86(99),\,44(105),\,59(105),\,
89(105),\,31(117),\\
&43(117),\,83(117),\,103(117),\,35(135),\,
43(135),\,88(135),\,26(143),\,86(143),\\
&125(143),\,35(165),\,
37(175),\,87(175),\,162(175),\,34(189),\,53(189),\,155(195),\\
&85(225),\,130(225),\,157(225),\,202(225),\,137(231),\,158(231),\,
104(273),\\
&146(273),\,188(273),\,65(275),\,175(275),\,152(297),\,
218(297),\,79(315),\\
&284(315),\,295(315),\,
87(325),\,112(325),\,162(325),\,16(351),\,44(351),\\
&97(351),\,286(351),\,313(351),\,15(385),\,
225(385),\,290(385),\,191(429),\\
&203(429),\,284(429),\,34(455),\,454(455),\,130(495),\,230(495),\,
395(495),\\
&179(525),\,362(525),\,445(525),\,494(525),\,
335(585),\,355(585),\,412(585),\\
&490(585),\,7(675),\,
232(675),\,277(675),\,502(675),\,200(693),\,257(693),\\
&515(693),\,445(715),\,500(715),\,555(715),\,356(819),\,538(819),\,
629(819),\\
&100(825),\,145(825),\,265(825),\,
475(825),\,179(945),\,494(945),\,562(975),\\
&637(975),\,662(975),\,
862(975),\,937(975),\,115(1001),\,808(1001),\,5(1155),\\
&809(1155),\,845(1155),\,950(1155),\,614(1287),\,742(1287),\,
1010(1287),\\
&767(1365),\,977(1365),\,1235(1365),\,
350(1485),\,220(1575),\,662(1575),\\
&1012(1575),\,1390(1575),\,470(1755),\,580(1755),\,610(1755),\,
880(1755),\\
&564(1925),\,949(1925),\,
1089(1925),\,1334(1925),\,1474(1925),\,1859(1925),\\
&202(2079),\, 895(2079),\,911(2079),\,1105(2145),\,
1670(2145),\,1012(2275),\\
&1362(2275),\,
1537(2275),\,647(2457),\,853(2457),\,1210(2457),\,1214(2457),\\
&2365(2457),\,2384(2457),\,670(2475),\,
2245(2475),\,2290(2475),\\
&2264(3003),\,1390(3465),\,
416(3861),\,3195(5005),\,1600(5775),\\
&2920(6435),\,7825(10395),\,583939(675675).
\endaligned$$

It is easy to check that the least common multiple of
$n_1,\ldots,n_{173}$ is the odd number
$$3^3\times 5^2\times7\times11\times 13=675675.$$
Since $A_1=\{a_s(n_s)\}_{s=1}^{173}$ covers $0,\ldots,675674$, it covers all the integers.

Using the software {\it Mathematica} and the main tables of
[BLSTW, pp.\,1--59], below we associate each
$n\in\{n_1,\ldots,n_{173}\}$ with $m_n$ distinct primitive prime
divisors $p_{n,1},\ldots,p_{n,m_n}$ of $2^n-1$ and write
$n:\,p_{n,1},\ldots,p_{n,m_n}$ for this, where $m_n$ is the number
of occurrences of $n$ among the moduli $n_1,\ldots,n_{173}$. For
those
$$n\in\{1485,\,3003,\,3465,\,3861,\,5005,\,5775,\,6435,\,10395,\,675675\},$$
as $m_n=1$ we just need one primitive prime divisor of $2^n-1$
whose existence is guaranteed by Bang's theorem; but they are too
large to be included in the following list.
\medskip

3: 7; \qquad 5: 31; \qquad 7: 127; \qquad 9: 73; \qquad 11: 23,
89;\qquad 13: 8191;

15: 151; \qquad\quad\ \, 21: 337; \qquad\qquad\ \ \, 25: 601,
1801;\qquad\  \ \,27: 262657;

33: 599479; \qquad 35: 71, 122921;\qquad 39: 79, 121369; \qquad
45: 631, 23311;

55: 881, 3191, 201961; \qquad 63: 92737, 649657;\qquad\qquad\ \
65: 145295143558111;

75: 100801, 10567201; \qquad 77: 581283643249112959;

91: 911, 112901153, 23140471537;\qquad 99: 199, 153649,
33057806959;

105: 29191, 106681, 152041;\qquad\qquad 117: 937, 6553, 86113,
7830118297;

135: 271, 348031, 49971617830801;

143: 724153, 158822951431, 5782172113400990737;

165: 2048568835297380486760231;

175: 39551, 60816001, 535347624791488552837151;

189: 1560007, 207617485544258392970753527;

195: 134304196845099262572814573351;

225: 115201, 617401, 1348206751, 13861369826299351;

231: 463, 4982397651178256151338302204762057;

273: 108749551, 4093204977277417, 86977595801949844993;

275: 382027665134363932751,
4074891477354886815033308087379995347151;

297: 8950393, 170886618823141738081830950807292771648313599433;

315: 870031, 983431, 29728307155963706810228435378401;

325: 7151, 51879585551,
46136793919369536104295905320141225322603397396- 44049093601;

351: 446473, 29121769, 571890896913727, 93715008807883087,
15083242680017- 3710177;

385: 55441, 1971764055031,
3105534168119044447812671975596513457115147- 3925765532041;

429: 17286204937, 1065107717756542892882802586807,
16783351554928582788- 5461382441449;

455: 200201,
477479745360834380327098898433221409835178252774757745602-
8391624903856636676854631;

495: 991, 334202934764737951438594746151,
60847771595376357965505368637- 41698483921;

525: 4201, 7351, 181165951,
32598550887552758766960709722266755711622113- 9090131514801;

585: 2400314671, 339175003117573351, 255375215316698521591,
272833453603- 4592865339299805712535332071;

675: 1605151, 1094270085398478390395590841401,
284249626318864764008979- 4561760551,
470390038503476855180627941942761032401;

693: 289511839, 2868251407519807,
3225949575089611556532995773813585269-
068981944367719218489696982054779837928902323497;

715: 249602191565465311, 598887853030285391,
40437156024702109576962112-
69051564057334878401893925719287086587582273263116838732848215441416415-

\noindent 0624064713711;

819: 2681001528674743, 219516331727145697249308031,
2149497317930391319-
0133458460563964459380529075838941297352657742148160962406273546512257;

825: 702948566745151, 9115784422509601,
4108316654247271397904922852177- 568560929751,
101249241260240615605217612230376981800142669401;

945: 124339521078546949914304521499392241,
893712833189249887135446424-
72309024678004403189516730060412595564942724011446583991926781827601;

975: 1951, 8837728285481551, 26155966684789722885001,
166376338119230863- 5718252801,
429450077043962550968970748284276205679121714346778186776993-
9979855730352201;

1001: 6007, 6952744694636960851412179090394909207;

1155: 2311, 6250631311, 494224324441,
2600788923312052743240883667728867-
90199621606534384599607578416912079166019131912393708208277038936454393-

\noindent 545946152508951;

1287: 216217, 71477407, 141968533929529744009;

1365: 469561, 52393016292934591, 2224981001722824694441;

1575: 82013401, 32758188751, 76641458269269601,
764384916291005220555242- 939647951;

1755: 3511, 196911, 4242734772486358591, 85488365519409100951;

1925: 11551, 13167001, 1891705201, 5591298184498951,
292615400703113951, 5627063397043739893603449551;

2079: 4159, 16633, 80932047967;

2145: 96001053721, 347878768688881;

2275: 218401, 28319200001, 1970116306308855665077103351;

2457: 565111, 1410319, 21287449, 41194063, 16751168775662428927,
178613107- 4995391292297656133027144291751;

2475: 4951, 143551, 1086033846151.

\medskip

Observe that $p_{n,j}>5$ for all $n\in\{n_1,\ldots,n_{173}\}$ and
$1\ls j\ls m_n$. In view of the above, Theorem 1.1 has been
proved. \qed

\heading{3. Proof of Theorem 1.2}\endheading

Recall that an odd prime $p$ is called a Wieferich prime if
$2^{p-1}\eq1\ (\mo\ p^2)$. The only known Wieferich primes are
1093 and 3511, and there are no others below $1.25\times 10^{15}$
(cf. [R, p.\,230]).

Suppose that $n\not=6$ is an integer greater than than one, and
$p$ is a primitive prime divisor of $2^n-1$. Then $n$ is the order
of $2$ mod $p$ and hence $p-1$ is a multiple of $n$ by Fermat's
little theorem. Thus $2^n-1\mid 2^{p-1}-1$, and hence $p^2\nmid
2^n-1$ if $p$ is not a Wieferich prime.

Let $A_1=\{a_s(n_s)\}_{s=1}^{173}$ and $p_1,\ldots,p_{173}$ be as
described in Theorem 1.1.
For each $s=1,\ldots,173$ let $q_s$ be a
primitive prime divisor of $2^{p_s^2}-1$. Then
$p_1,\ldots,p_{173},$ $q_1,\ldots,q_{173}$ are distinct odd primes
since $\{p_1^2,\ldots,p_{173}^2\}\cap\{n_1,\ldots,n_{173}\}=\em$.

For each $s=1,\ldots,173$ let $\al_s$ be the largest positive
integer with $p_s^{\al_s}\mid 2^{n_s}-1$. Since $3511$ is the only
Wieferich prime in the set $\{p_1,\ldots,p_{173}\}$, we have
$\al_s=1$ if $p_s\not=3511$. In the case $p_s=3511$, we have
$\al_s=2$ since $3511^2\mid 2^{3510}-1$ but $3511^3\nmid
2^{3510}-1$.

Let $M=2^{2L}\prod_{s=1}^{173}p_s^{\al_s+2}q_s$, where  $L$ is the
smallest positive integer satisfying
$$2^L-1>\max\{16N,p_1^{\al_1+1},\ldots,p_{173}^{\al_{173}+1}\}.$$
By the Chinese Remainder Theorem, there exists a unique
$a\in\{1,\ldots,M\}$ such that
$$1+3\cdot2^L(2^{2L})\cap\bigcap_{s=1}^{173}\l(x_s^{b_s}(p_s^{\al_s+2})\cap y_s^{b_s}(q_s)\r)=a(M).$$

Let $m\ls N$ be a positive integer relatively prime to
$3\cdot5\cdot7\cdot11\cdot 13=15015$, and write $m=2^\al m_0$ with
$\al\in\N$, $m_0\in\Z^+$ and $2\nmid m_0$. Let
$s\in\{1,\ldots,173\}$. Since $n_s$ is a divisor of
$3^3\cdot5^2\cdot7\cdot11\cdot 13=675675$, we have $\gcd(m,n_s)=1$
and hence $m_0b_s\eq a_s\ (\mo\ n_s)$ for some $b_s\in\N$.

As the order of $2$ mod $p_s$ is the odd number $n_s$,
$n_s$ divides $(p_s-1)/\gcd(2^\al,p_s-1)$
and hence
$$2^{(p_s-1)/\gcd(2^\al,\ p_s-1)}\eq 1\ (\mo\ p_s),
\, 2^{p_s(p_s-1)/\gcd(2^\al,\ p_s-1)}\eq 1\ (\mo\ p_s^2),
\,\ldots.$$
Since there is a primitive root modulo $p_s^{\al_s+2}$ and
$$2^{\varphi(p_s^{\al_s+2})/\gcd(2^\al,\varphi(p_s^{\al_s+2}))}
=2^{p_s^{\al_s+1}(p_s-1)/\gcd(2^\al,p_s-1)}\eq 1\ (\mo\
p_s^{\al_s+2})$$ (where $\varphi$ is Euler's totient function), by
[IR, Proposition 4.2.1] there exists $x_s\in\Z$ with
$x_s^{2^\al}\eq 2\ (\mo\ p_s^{\al_s+2})$. Similarly, the order
$p_s^2$ of $2$ mod $q_s$ divides $(q_s-1)/\gcd(2^\al,q_s-1)$,
therefore $2^{(q_s-1)/\gcd(2^\al,q_s-1)}\eq1\ (\mo\ q_s)$ and
hence $y_s^{2^\al}\eq 2\ (\mo\ q_s)$ for some $y_s\in\Z$.

Let $x\gs0$ be an element of $a(M)$. As $A_1$ is a cover of $\Z$,
for any $n\in\N$ there is an $s\in\{1,\ldots,173\}$ such that
$n\eq a_s\ (\mo\ n_s)$. Clearly
$$x^m\eq (x_s^{b_s})^m=(x_s^{2^\al})^{m_0b_s}\eq2^{m_0b_s}\ (\mo\ p_s^{\al_s+2}),$$
thus
$$x^m-2^n\eq 2^{m_0b_s}-2^{a_s}\eq0\ (\mo\ p_s^{\al_s})$$
since $2^{n_s}\eq1\ (\mo\ p_s^{\al_s})$ and $m_0b_s\eq a_s\ (\mo\ n_s)$.

As $16m\ls 16N<2^L-1$ and $x\eq1+3\cdot2^L\ (\mo\ 2^{2L})$, we have $|x^m-2^n|\gs 2^L-1>p_s^{\al_s+1}$
by [C, Lemma 1]. So $|x^m-2^n|\not=0,p_s^{\al_s},p_s^{\al_s+1}$.
If $x^m-2^n$ is not divisible by $p_s^{\al_s+2}$, then it must
have at least two distinct prime divisors.

Now we assume that $x^m-2^n\eq0\ (\mo\ p_s^{\al_s+2})$. Note that
$2^n\eq x^m\eq 2^{m_0b_s}\ (\mo\ p_s^{\al_s+2})$.
Since $n_s$ is the order of $2$ mod $p_s^{\al_s}$
and not the order of $2$ mod $p_s^{\al_s+1}$, by [C, Corollary 3]
we have $2^n\eq 2^{m_0b_s}\ (\mo\ q_s)$. Thus
$$x^m-2^n\eq(y_s^{b_s})^{2^\al m_0}-2^{m_0b_s}\eq0\ (\mo\ q_s)$$
and so the nonzero integer $x^m-2^n$ has at least two distinct
prime divisors (including $p_s$ and $q_s$).

By the above, we have proved the desired result. \qed

\medskip

\Remark\ 3.1. Given $m,n\in\Z^+$ and an odd prime $p$, the
equation $x^m-2^n=p^b$ with $b,x\in\N$ only has finitely many
solutions. As observed by the referee, this is a consequence of
the Darmon-Granville theorem in [DG]. In the case $m=2$, all the
finitely many solutions are effectively computable by the
algorithms given by Weger [W].

\medskip

\heading{4. Proof of Theorem 1.3}\endheading

\proclaim{Lemma 4.1} Let $c\in\Z^+$, and define $\{U_n\}_{n\gs0}$ by
$$U_0=0,\ U_1=1,\ \t{and}\ U_{n+1}=cU_{n}+U_{n-1}\ \t{for}\ n=1,2,3,\ldots.$$
Suppose that $n>0$ is an integer with $n\eq2\ (\mo\ 4)$ and $p$
is a prime divisor of $U_n$ which
divides none of $U_1,\ldots,U_{n-1}$. Then $U_{kn+r}\equiv U_r\ (\mo\ p)$ for all
$k\in\N$ and $r\in\{0,\ldots,n-1\}$.
\endproclaim

\Proof. By [HS, Lemma 2],
$U_{n+1}\equiv -(-1)^{n/2}=1\ (\mo\ p)$. If $k\in\N$ and $r\in\{0,\ldots,n-1\}$,
then $U_{kn+r}\equiv U^k_{n+1}U_r\ (\mo\ U_n)$ by [HS, Lemma 3] or [S92, Lemma 2],
therefore $U_{kn+r}\eq U_r\ (\mo\ p)$. \qed

\medskip
\noindent{\it Proof of Theorem 1.3}.
Let $b_1(m_1),\ldots,b_{24}(m_{24})$ be the following 24 residue classes:
$$\aligned
&1(3),\,2(5),\,3(5),\,4(7),\,6(7),\,0(9),\,5(15),\,
11(15),\,9(21),\,12(21),\\
&1(35),\,14(35),\,24(35),\,29(35),\,6(45),\,15(45),\,
29(45),\,30(45),\\
&5(63),\,23(63),\,44(63),\,66(105),\,21(315),\,89(315).
\endaligned$$
It is easy to check that $\{b_t(m_t)\}_{t=1}^{24}$ forms a cover
of $\Z$ with odd moduli. Set $m_0=1$. Then
$$B=\{1(2m_0), 2b_1(2m_1),\ldots,2b_{24}(2m_{24})\}$$
is a cover of $\Z$ with all the moduli congruent to 2 mod 4.

Let $u_n=F_{3n}/2$ for $n\in\N$.
As we mentioned in Section 1, $u_0=0$, $u_1=1$ and
$u_{n+1}=4u_n+u_{n-1}$ for $n=1,2,3,\ldots$.
For a prime $p$ and an integer $n>0$, we call $p$
a {\it primitive prime divisor of $u_n$} if $p\mid u_n$ but $p\nmid u_k$ for those $0<k<n$.

Let $p_0,\ldots,p_{24}$ be the following 25 distinct primes respectively:
$$\aligned
&2,\,19,\,31,\,11,\,211,\,29,\,5779,\,541,\,
181,\,31249,\,1009,\,767131,\,21211,\,911,\\
&71,\,119611,\,42391,\,271,\,811,\,379,\,
912871,\,85429,\,631,\,69931,\,17011.
\endaligned$$
One can easily verify that each $p_t\ (0\ls t\ls 24)$
is a primitive prime divisor of $u_{2m_t}$.

The residue class $a(M)$ in Theorem 1.3 is actually the intersection
of the following 25 residue classes with the moduli
$p_0,\ldots,p_{24}$ respectively:
$$\aligned
&1(2),\ 2(19),\ 14(31),\ 4(11),\ 94(211),\
5(29),\  0(5779),\ 156(541),\ 76(181),
\\&10727(31249),\ 501(1009),\ 2(767131),\ 7199(21211),\ 257(911),\ 30(71),
\\&13909(119611),\ 9054(42391),\ 85(271),\ 292(811),\
72(379),\  80065(912871),
\\&40368(85429),\
205(631),\ 19928(69931),\ 497(17011).
\endaligned$$

It is known that the only solutions of the diophantine equation
$F_n=2x^2$ with $n,x\in\N$ are $(n,x)=(0,0),(3,1),(6,2)$. (Cf.
[Co, Theorem 4].) Let $x$ be any integer in the residue class
$a(M)$. Then $|x|>2$ and hence $x^2\not=u_n=F_{3n}/2$ for all
$n\in\N$. With the help of Lemma 4.1 in the case $c=4$, one can
check that $x^2\eq u_1=1\ (\mo\ p_0)$ and $x^2\eq u_{2b_t}\ (\mo\
p_t)$ for all $t=1,\ldots,24$.

Let $n$ be any nonnegative integer. As $B$ forms a cover of $\Z$,
$n\eq 1\ (\mo\ 2m_0)$ or $n\eq 2b_t\ (\mo\ 2m_t)$ for some $1\ls
t\ls 24$. By Lemma 4.1 with $c=4$, if $n\eq1\ (\mo\ 2m_0)$ then
$u_n\eq u_1=1\ (\mo\ p_0)$ and hence $x^2-u_n\eq x^2-1\eq0\ (\mo\
p_0)$; if $n\eq 2b_t\ (\mo\ 2m_t)$ then $u_n\eq u_{2b_t}\ (\mo\
p_t)$ and hence $x^2-u_n\eq x^2-u_{2b_t}\eq0\ (\mo\ p_t)$. Thus,
it remains to show that for any given $a,b\in\N$ we can deduce a
contradiction if $x^2-u_{1+2m_0a}=\pm 2^b$ or
$x^2-u_{2b_t+2m_ta}=\pm p_t^b$ for some $1\ls t\ls 24$.
\medskip

{\tt Case 4.0}. $x^2-u_{1+2a}=\pm 2^b$.

As $p_2=31$ and $p_3=11$ are primitive prime divisors of $u_{2m_2}=u_{2m_3}=u_{10}$,
and
$$u_1=1,\ u_3=17,\ u_5=305,\ u_7=5473,\ u_9=98209$$
have residues $1,-14,-5,-14,1$ modulo $31$ and residues
$1,-5,-3,-5,1$ modulo $11$ respectively. If $2a+1\not\eq 5\ (\mo\
10)$, then by Lemma 4.1 we have
$$x^2-u_{1+2a}\eq 10-1,10-(-14)\not\eq \pm 1,\pm 2,\pm 4,\pm 8,\pm 16\ (\mo\ 31)$$
which contradicts $x^2-u_{1+2a}=\pm 2^b$. (Note that $2^5\eq 1\
(\mo\ 31)$.) So $2a+1\eq 5\ (\mo\ 10)$. It follows that
$$x^2-u_{1+2a}\eq 10-(-5)\eq -2^4\ (\mo\ 31)\ \text{and}\
x^2-u_{1+2a}\eq 5-(-3)=2^3\ (\mo\ 11).$$ Thus $x^2-u_{1+2a}$ can
only be $-2^b$ with $b\eq 4\ (\mo\ 5)$, which cannot be congruent to
$2^3$ mod 11. (Note that $2^5\eq -1\ (\mo\ 11)$.) So we have a
contradiction.
\medskip

{\tt Case 4.1}. $x^2-u_{2+6a}=\pm 19^b$.

Observe that
$$u_0=0,\ u_2=4,\ u_4=72,\ u_6=1292,\ u_8=23184$$
have residues $0,4,-5,5,-4$ modulo $11$ and $0,4,10,-10,-4$ modulo
$31$ respectively. Also, $19^b\eq 2^{3b}\eq\pm 1,\pm 2,\pm 4,\pm
3,\pm 5\ (\mo\ 11)$ and $19^5\eq(-2^2\cdot3)^5\eq-3^5\eq5\ (\mo\ 31)$.

If $2+6a\eq 0\ (\mo\ 10)$, then
$$x^2-u_{2+6a}\eq 5-0\eq 19^8,-19^3\ (\mo\ 11)$$
and hence
$x^2-u_{2+6a}=(-1)^{d-1}19^{3+5d}$ for some $d\in\N$,
this leads to a contradiction since $x^2-u_{2+6a}\eq 10-0\ (\mo\ 31)$
but
$$19^{3+5d}\eq 8\times 5^d\eq 8,9,14\not\eq\pm 10\ (\mo\ 31).$$

Now we handle the case $2+6a\eq 2\ (\mo\ 10)$.
Since 181 is a primitive prime divisor of $u_{30}$, and $6a\eq0\ (\mo\ 30)$ and $19^2\eq-1\ (\mo\ 181)$,
we have
$$x^2-u_{2+6a}\eq76^2-u_2\eq-20\not\eq\pm 19^b\ (\mo\ 181)$$
which leads a contradiction.

 If $2+6a\eq 4\ (\mo\ 10)$, then
 $x^2-u_{2+6a}\eq 10-10=0\ (\mo\ 31)$.
If $2+6a\eq 6\ (\mo\ 10)$, then $x^2-u_{2+6a}\eq 5-5=0\ (\mo\
11)$. So, when  $2+6a\eq 4,6\ (\mo\ 10)$ we get a contradiction
since $x^2-u_{2+6a}=\pm 19^b$.

If $2+6a\eq 8\ (\mo\ 10)$, then
$x^2-u_{2+6a}\eq 5-(-4)\eq 19^2,-19^7\ (\mo\ 11)$
and hence $x^2-u_{2+6a}=(-1)^d19^{2+5d}$ for some $d\in\N$,
this leads a contradiction since $x^2-u_{2+6a}\eq 10-(-4)\eq-11\times10\ (\mo\
31)$ but
$$19^{2+5d}\eq-11\times 5^d\eq-11,-11\times 5,-11\times(-6)\not\eq\pm 11\times 10\ (\mo\ 31).$$
\medskip

{\tt Case 4.2}. $x^2-u_{4+10a}=\pm 31^b$.

As $x^2-u_{4+10a}\eq 5-u_4\eq 5-(-5)\eq -1\ (\mo\ 11)$ and $31^b\eq
(-2)^b\eq 1,-2,4,-8,16\ (\mo\ 11)$, we must have
$x^2-u_{4+10a}=-31^b$ with $b\eq 0\ (\mo\ 5)$. As $31^5\eq 2^3=8\
(\mo\ 19)$, $31^b\eq 8^{b/5}\eq\pm 1,\pm 8,\pm 7\ (\mo\ 19)$. If
$3\nmid a$, then $4+10a\eq 0,2\ (\mo\ 6)$ and hence
$$x^2-u_{4+10a}\eq 4-u_0,4-u_2\eq 4,0\not\eq -31^b\ (\mo\ 19).$$
Thus $a=3c$ for some $c\in\N$. As
$$-8^{b/5}\eq -31^b=x^2-u_{4+10a}\eq 4-u_4=4-72\eq 8\ (\mo\ 19),$$
we have $b/5-1\eq 3\ (\mo\ 6)$ and hence $b=20+30d$ for some
$d\in\N$. As $31^{10}\eq -1\ (\mo\ 181)$, we have
$31^b=31^{20+30d}\eq (-1)^{2+3d}=(-1)^d\ (\mo\ 181)$. On the other
hand,
$$-31^b=x^2-u_{4+10a}=x^2-u_{4+30c}\eq 76^2-u_4\eq -16-72=-88\ (\mo\ 181).$$
So we get a contradiction.
\medskip

{\tt Case 4.3}. $x^2-u_{6+10a}=\pm 11^b$.

As $x^2-u_{6+10a}\eq 10-(-10)\eq -11\ (\mo\ 31)$, and the order of
11 mod 31 is 30, we have $x^2-u_{6+10a}=(-1)^{d-1}11^{1+15d}$ for some $d\in\N$. Since
$11^{15}\eq (-8)^{15}=(-2^9)^5\eq 1\ (\mo\ 19)$, we have
$x^2-u_{6+10a}\eq\pm 11\ (\mo\ 19)$.

If $6+10a\eq 0,2\ (\mo\ 6)$, then
$$x^2-u_{6+10a}\eq 4-u_0,4-u_2\not\eq\pm 11\ (\mo\ 19).$$
So $6+10a\eq 4\ (\mo\ 6)$, i.e., $a=1+3c$ for some $c\in\N$.
Therefore
$$x^2-u_{6+10a}=x^2-u_{16+30c}\eq -16-u_{16}\eq -16-47\eq -11\times88\ (\mo\ 181).$$
Note that
$$(-11)^{15d}\eq (-49)^d\eq 1,\,-49,\,48\not\eq88\ (\mo\ 181).$$
As $x^2-u_{6+10a}=(-11)^{1+15d}$, we get a contradiction.
\medskip

{\tt Case 4.4}. $x^2-u_{8+14a}=\pm 211^b$.

As $p_5=29$ is a primitive divisor of $u_{2m_5}=u_{14}$, we have
$x^2-u_{8+14a}\eq 25-u_8\eq 25-13=12\ (\mo\ 29)$.

Since $2$ is a primitive root mod $29$, $211\eq 2^3\ (\mo\ 29)$,
$2^{3\times21}\eq 2^7\eq 12\ (\mo\ 29)$, and $2^{3\times 7}\eq
12^3\eq-12\ (\mo\ 29)$, we have $x^2-u_{8+14a}=(-1)^{d-1}211^{7+14d}$ for some $d\in\N$.

Observe that
$$\aligned
x^2-u_{8+14a}&\eq 10-u_0,10-u_2,10-u_4,10-u_6,10-u_8,\\
&\eq 10-0,10-4,10-10,10-(-10),10-(-4)\ (\mo\ 31).
\endaligned$$
Clearly $211\eq 5^2\ (\mo\ 31)$ and $5^3\eq1\ (\mo\ 31)$, thus
$$211^{7+14d}\eq 5^{14+28d}\eq 5^{2+d}\eq-6,\,1,\,5\ (\mo\ 31).$$
Therefore $2\mid d$, $3\mid d$ and $8+14a\eq 2\ (\mo\ 10)$.
It follows that $a=1+5c$ for some $c\in\N$ and $d=6e$ for some $e\in\N$.

Note that $$x^2-u_{8+14(1+5c)}\eq x^2-u_2\eq 5-4=1\ (\mo\ 11)$$
and
$$(-1)^{d-1}211^{7+14d}\eq -2^{7(1+12e)}\eq -2^{7(1+2e)}\ (\mo\ 11).$$
So $2^{7(1+2e)}\eq-1\eq 2^5\ (\mo\ 11)$,
hence $7(1+2e)\eq5\eq35\ (\mo\ 10)$ and thus $e\eq2\ (\mo\ 5)$. Therefore
$7+14d\eq 7+84\times2\eq35\ (\mo\ 140)$ and hence
$$211^{7+14d}\eq(-2)^{35}\eq\l(\f{-2}{71}\r)=-\l(\f 2{71}\r)=-1\ (\mo\ 71)$$
by the theory of quadratic residues, but
$$x^2-u_{8+14a}=x^2-u_{22+70c}\eq 30^2-u_{22}=900-13888945017644\eq14\ (\mo\ 71),$$
so we get a contradiction from the equality $x^2-u_{8+14a}=-211^{7+14d}$.
\medskip

{\tt Case 4.5}. $x^2-u_{12+14a}=\pm 29^b$.

As $29^b\eq (-2)^b\eq \pm 1,\pm 2,\pm4,\pm 8,\pm 16\ (\mo\ 31)$,
$x^2\eq 14^2\eq10\ (\mo\ 31)$ and
$$u_{12+14a}\eq u_0,u_2,u_4,u_6,u_8\eq 0,4,10,-10,-4\ (\mo\ 31),$$
we have $x^2-u_{12+14a}\not\eq\pm29^b\ (\mo\ 31)$. So, a
contradiction occurs.
\medskip

{\tt Case 4.6}. $x^2-u_{0+18a}=\pm 5779^b$.

As $x^2-u_{18a}\eq 2^2-u_0=4\ (\mo\ 19)$, $5779\eq 3\ (\mo\ 19)$
and the order of $3$ mod $19$ equals 18,
we have $x^2-u_{18a}=(-1)^{d-1}5779^{5+9d}=(-5779)^{5+9d}$ for some $d\in\N$.

Note that
$$(-5779)^{9d}\eq (-13)^{9d}\eq(2^2)^{3d}\eq 2^d\eq 1,\,2,\,4,\,8,\,16\ (\mo\ 31)$$
and $5779^5\eq13^5\eq 6\ (\mo\ 31)$.
Thus
$$x^2-(-5779)^{5+9d}\eq 10+6\times 2^d\eq -15,\,-9,\,3,\,-4,\,13\ (\mo\ 31)$$
while $u_{18a}\eq u_0,u_2,u_4,u_6,u_8\eq 0,4,10,-10,-4\ (\mo\ 31)$.
As $u_{18a}=x^2-(-5779)^{5+9d}$, we must have
$18a\eq 8\ (\mo\ 10)$ and $d=3+5e$ for some $e\in\N$.

Observe that $x^2-u_{18a}\eq5-u_8\eq-2\ (\mo\ 11)$ but
$$(-5779)^{5+9d}\eq(-2^2)^{5+9(3+5e)}=(-1)^e2^{64+90e}\eq(-1)^e2^4\not\eq-2\ (\mo\ 11).$$
So a contradiction occurs.
\medskip

{\tt Case 4.7}. $x^2-u_{10+30a}=\pm 541^b$.

As $x^2-u_{10+30a}\eq 5-u_0\eq 5\ (\mo\ 11)$ and
$$541^b\eq2^b\eq \pm 1,\pm 2,\pm 3,\pm 4,\pm 8,\pm 16\ (\mo\ 11),$$
$x^2-u_{10+30a}=(-1)^d541^{4+5d}$ for some $d\in\N$,
and hence we have a contradiction since
 $x^2-u_{10+30a}\eq 10-u_0=10\ (\mo\ 31)$ but
$$541^{4+5d}\eq (2\times7)^{4+5d}\eq7\times 5^d\eq7,\,7\times5,\,7\times(-6)\not\eq\pm10\ (\mo\ 31).$$
\medskip

{\tt Case 4.8}. $x^2-u_{22+30a}=\pm 181^b.$

As $x^2-u_{22+30a}\eq 5-u_2=5-4\ (\mo\ 11)$ and
$181^b\eq 5^b\eq 1,5,3,4,-2\ (\mo\ 11)$,
we have $x^2-u_{22+30a}=181^b$ with $b=5d$ for some $d\in\N$.
Since
$x^2-u_{22+30a}\eq x^2-u_2\eq 10-4=6\ (\mo\ 31)$ and $181^{5d}\eq
(-5)^{5d}\eq6^d\eq 1,6,5,-1,-6,-5\ (\mo\ 31)$, $d=1+6e$ for some $e\in\N$.
Note that $x^2-u_{22+30a}\eq 4-u_4=4-72\eq 8\ (\mo\ 19)$ but
$$181^{5d}\eq(-9)^{5d}=(-3^{10})^d\eq 3^d
=3^{1+6e}\eq 3\times7^e\eq 3,\,2,\,-5\not\eq 8\ (\mo\ 19).$$
\medskip

{\tt Case 4.9}. $x^2-u_{18+42a}=\pm 31249^b.$

Note that $31249^b\eq 1^b=1\ (\mo\ 31)$, $x^2\eq 10\ (\mo\ 31)$ and also
$$u_{18+42a}\eq u_0,u_2,u_4,u_6,u_8\eq0,\,4,\,10,\,-10,\,-4\ (\mo\ 31).$$
Therefore $x^2-u_{18+42a}\not\eq\pm31249^b\ (\mo\ 31)$.
\medskip

{\tt Case 4.10}. $x^2-u_{24+42a}=\pm 1009^b.$

As $x^2-u_{24+42a}\eq 4-u_0=4\ (\mo\ 19)$, $1009\eq 2\ (\mo\ 19)$
and $2$ is a primitive root modulo $19$,
we have $x^2-u_{24+42a}=(-1)^d1009^{2+9d}=(-1009)^{2+9d}$ for some $d\in\N$.
Observe that $u_{10}=416020$ and $x^2-u_{24+42a}\eq 25-u_{10}\eq 10\ (\mo\ 29)$.
But $6^7\eq-1\ (\mo\ 29)$ and hence
$$(-1009)^{2+9d}\eq6^{2+9d}\eq \pm1,\,\pm6,\,\pm7,\,\pm13,\,\pm9,\,\pm4,\,\pm5\not\eq 10\ (\mo\ 29).$$
So we get a contradiction.
\medskip

{\tt Case 4.11}. $x^2-u_{2+70a}=\pm 767131^b$.

Observe that $x^2-u_{2+70a}\eq 5^2-u_2\eq-8\ (\mo\ 29)$ and
$$767131^b\eq(-6)^b\eq1,\,-6,\,7,\,-13,\,-9,\,-4,\,-5\ (\mo\ 29).$$
So a contradiction occurs.
\medskip

{\tt Case 4.12}. $x^2-u_{28+70a}=\pm 21211^b$.

As $28+70a\eq0\ (\mo\ 14)$, we have $x^2-u_{28+70a}\eq 5^2-u_0\eq-4\ (\mo\ 29)$.
On the other hand,
$$21211^b\eq \pm 12^b\eq \pm1,\,\pm12\ (\mo\ 29).$$
Thus we have a contradiction.
\medskip

{\tt Case 4.13}. $x^2-u_{48+70a}=\pm911^b$.

Note that $x^2-u_{48+70a}\eq 5^2-u_6=25-1292\eq 9\ (\mo\ 29)$
but $911^b\eq 12^b\eq \pm1,\,\pm12\ (\mo\ 29)$.
\medskip

{\tt Case 4.14}. $x^2-u_{58+70a}=\pm 71^b$.

Observe that $x^2-u_{58+70a}\eq 5^2-u_2\eq-8\ (\mo\ 29)$ but
$$71^b\eq 13^b\eq\pm 1,\,\pm 13,\,\pm5,\,\pm7,\,\pm 4,\,\pm 6,\,\pm 9\ (\mo\ 29).$$
\medskip

{\tt Case 4.15}. $x^2-u_{12+90a}=\pm 119611^b.$

Since $x^2-u_{12+90a}\eq 4-u_0\eq 4\ (\mo\ 19)$ and
$$119611^b\eq 6^b\eq1,\,6,\,-2,\,7,\,4,\,5,\,-8,\,9,\,-3\ (\mo\ 19),$$
we must have
$x^2-u_{12+90a}=119611^b$ with $b=4+9d$ for some $d\in\N$.
Note that $x^2-u_{12+90a}\eq10-u_2=6\eq10\times13\ (\mo\ 31)$, but
$$119611^{4+9d}\eq 13^{4+9d}\eq 10(-2)^d\ (\mo\ 31)$$
with $(-2)^d\eq\pm1,\,\pm2,\,\pm4,\,\pm8,\,\pm16\not\eq13\ (\mo\ 31)$. So
we have a contradiction.
\medskip

{\tt Case 4.16}. $x^2-u_{30+90a}=\pm 42391^b.$

As $x^2-u_{30+90a}\eq 4-u_0=4\ (\mo\ 19)$, $42391\eq 2\ (\mo\
19)$ and $2$ is a primitive root mod 19,
we have $x^2-u_{30+90a}=(-1)^d42391^{2+9d}$ for some $d\in\N$.

 Note that $x^2-u_{30+90a}\eq 10-0\ (\mo\ 31)$ and
$$(-42391)^{2+9d}\eq(-14)^{2+9d}\eq10(-2^4)^{3d}\eq10(-1)^d2^{2d}\ (\mo\ 31).$$
Since the only residues of powers of 2 modulo 31 are $1,\,2,\,4,\,8,\,16$, we must have
$x^2-u_{30+90a}=(-42391)^{2+9d}$ with $d$ divisible by both $5$ and $2$.
Write $d=10e$ with $e\in\N$. Then
$$x^2-u_{30+90a}=42391^{2+90e}\eq (-3)^{2+90e}\eq9\ (\mo\ 11),$$
which contradicts the fact $x^2-u_{30+90a}\eq5-u_0=5\ (\mo\ 11)$.
\medskip

{\tt Case 4.17}. $x^2-u_{58+90a}=\pm271^b$.

 Note that $x^2-u_{58+90a}\eq 10-u_8\eq14\ (\mo\ 31)$ while
 $$271^b\eq (-2)^{3b}\eq \pm 1,\,\pm 2,\,\pm 4,\,\pm 8,\,\pm 16\ (\mo\ 31).$$
\medskip

{\tt Case 4.18}. $x^2-u_{60+90a}=\pm 811^b$.

As $x^2-u_{60+90a}\eq 10-u_0=10\ (\mo\ 31)$ and
$811^b\eq5^b\eq 1,\,5,\,25\ (\mo\ 31)$.
we have a contradiction.
\medskip

{\tt Case 4.19}. $x^2-u_{10+126a}=\pm 379^b$.

Note that $x^2-u_{10+126a}\eq 2^2-u_4=4-72\eq8\ (\mo\ 19)$ but
$379^b\eq(-1)^b\eq\pm1\ (\mo\ 19)$.
\medskip

{\tt Case 4.20}. $x^2-u_{46+126a}=\pm 912871^b$.

Since $x^2-u_{46+126a}\eq2^2-u_4\eq 2^3\ (\mo\ 19)$,
$912871^b\eq2^{4b}\ (\mo\ 19)$ and the order of $2$ mod 19 is 18,
we must have $x^2-u_{46+126a}=-912871^{b}$ with $b=3+9d$ for some $d\in\N$.
Note that $x^2-u_{46+126a}\eq 5^2-u_4=25-72\eq11\ (\mo\ 29)$ but
$$912871^{3+9d}\eq 3^{2(3+9d)}\eq4^{1+3d}\eq
\pm1,\,\pm4,\,\pm13,\,\pm6,\,\pm5,\,\pm9,\,\pm7\ (\mo\ 29).$$
So we have a contradiction.
\medskip

{\tt Case 4.21}. $x^2-u_{88+126a}=\pm 85429^b.$

Observe that $x^2-u_{88+126a}\eq 5^2-u_4\eq11\ (\mo\ 29)$ but
$$85429^{b}\eq (-5)^b\eq1,\,-5,\,-4,\,-9,\,-13,\,7,\,-6\ (\mo\ 29).$$
So a contradiction occurs.
\medskip

{\tt Case 4.22}. $x^2-u_{132+210a}=\pm 631^b.$

Note that $x^2-u_{132+210a}\eq 4^2-u_2\eq1\ (\mo\ 11)$ and $631\eq2^2\ (\mo\ 11)$.
Since $2^5\eq-1\ (\mo\ 11)$ and $2^{10}\eq1\ (\mo\ 11)$, we must have
$x^2-u_{132+210a}=631^b$ with $b=5d$ for some $d\in\N$.
As $x^2-u_{132+210a}\eq 10 -u_2=6\ (\mo\ 31)$,
$631^5\eq (-2^2\times5)^5\eq-5^2\eq6\ (\mo\ 31)$
and the order of 6 mod 31 is 6, we can write $d=1+6e$ with $e\in\N$.
Thus
$$x^2-u_{132+210a}=631^{5+30e}\eq (2^2)^{5+30e}\eq (-2)^{1+6e}\eq-2,\,5,\,-3\ (\mo\ 19).$$
On the other hand, $x^2-u_{132+210a}\eq 4-u_0=4\ (\mo\ 19)$. This
leads to a contradiction.
\medskip

{\tt Case 4.23}. $x^2-u_{42+630a}=\pm 69931^b$.

As $42+630a\eq0\ (\mo\ 6)$, we have $x^2-u_{42+630a}\eq 2^2-u_0=4\ (\mo\ 19)$.
On the other hand, $69931^b\eq (-2)^{3b}\eq 1,\,-8,\,7\ (\mo\ 19)$.
So we get a contradiction.
\medskip

{\tt Case 4.24}. $x^2-u_{178+630a}=\pm 17011^b.$

Since $178+630a\eq10\ (\mo\ 14)$, we have
$$x^2-u_{178+630a}\eq 5^2-u_{10}=25-416020\eq10\ (\mo\ 29).$$
Note that $17011^b\eq(-12)^b\eq \pm1,\,\pm12\ (\mo\ 29)$.
So a contradiction occurs.

\medskip

In view of the above, we have completed the proof of Theorem 1.3. \qed

\medskip

\Ack. The authors would like to thank the referee for some helpful
comments.

\enddocument